\documentclass[preprint,sort&compress,final,12pt]{elsarticle}

\usepackage{amsmath}
\usepackage{amssymb}
\usepackage{graphicx}
\usepackage{latexsym}
\usepackage{epstopdf}
\usepackage{natbib}
\usepackage{color}
\usepackage{hyperref}

\newtheorem{theorem}{Theorem}[section]

\newtheorem{lemma}[theorem]{Lemma}

\numberwithin{equation}{section}

\def\Proof{\noindent{\bf Proof.}~}
\def\qed{\hfill$\square$\smallskip}
\def\dint{\displaystyle\int}

\def\dsup{\displaystyle\sup}

\def\dsum{\displaystyle\sum}
\def\dint{\displaystyle\int}
\def\dsup{\displaystyle\sup}

\def\dfrac#1#2{\frac{\displaystyle {#1}}{\displaystyle {#2}}}

\journal{\empty}
\date{}

\begin{document}

\begin{frontmatter}

\title{The coexistence of quasi-periodic and blow-up solutions in a superlinear Duffing equation\footnote{ Partially supported by the NSFC (11571041) and  the Fundamental Research Funds for the Central Universities.}}

\author[au1,au3]{Yanmei  Sun}

\ead[au1]{sunyanmei2009@126.com}

\author[au1]{Xiong Li\footnote{Corresponding author.}}

\address[au1]{School of Mathematical Sciences, Beijing Normal University, Beijing 100875, P.R. China.}

\address[au3]{School of Mathematics and Information Sciences, Weifang University, Weifang, Shandong, 261061, P.R. China.}

\ead[au1]{xli@bnu.edu.cn}

\begin{abstract}

In this paper we will construct a continuous positive periodic function $p(t)$ such that the corresponding superlinear Duffing equation
$$
x''+a(x)\,x^{2n+1}+p(t)\,x^{2m+1}=0,\ \ \ \ n+2\leq 2m+1<2n+1
$$
possesses a solution which escapes to infinity in some finite time, and also has infinitely many subharmonic and quasi-periodic solutions, where the coefficient $a(x)$ is an arbitrary positive smooth  periodic function defined in the whole real axis.
\end{abstract}


\end{frontmatter}

\section{Introduction}

In the early 1960's, Littlewood \cite{Littlewood} asked whether all solutions of  the second order differential equation
\begin{equation}\label{V}
x''+V_x(x,t)=0
\end{equation}
are bounded for all time, that is, $\dsup_{t\in \mathbb{R}} (|x(t)| + |x'(t)|) < +\infty$ holds for all
solutions $x(t)$ of Eq.(\ref{V}).

For the Littlewood boundedness problem, during the past years, people have
paid more attention to the  following equation with  the polynomial potentials
\begin{equation}\label{Duffing}
x''+x^{2n+1}+\dsum_{k=0}^{2n}p_k(t)x^k=0,
\end{equation}
where $p_k(t+1)=p_k(t)\ (k=0,1,\cdots, 2n)$,  since
$$
x'' + x^{2n+1} = 0
$$
is a very nice time-independent integrable system, of which all solutions are
periodic. Thus if $|x|$ is large enough, Eq.(\ref{Duffing}) can be treated as a perturbation
of an integrable system, then Moser's twist theorem could be applied to prove
the boundedness of all solutions.

The first result was due to Morris \cite{Morris}, who proved that all solutions of the equation with the
biquadratic potential
$$
x'' + 2x^3 = p(t) = p(t + 1)
$$
are bounded.

Using the famous Moser's twist theorem \cite{ Moser}, Diekerhoff and Zehnder \cite{Diekerhoff}
generalized Morris's results to Eq.(\ref{Duffing}). In  \cite{Diekerhoff}, the coefficients $p_k(t)$ are
required to be sufficiently smooth to construct a series  of variable changes to
transform Eq.(\ref{Duffing}) into a nearly integrable systems for large energies. In fact,
in \cite{Diekerhoff}, the smoothness on $p_k(t)$ depends on the index $k$.

Later, Yuan \cite{Yuan} proved that all solutions of  Eq.(\ref{Duffing}) are bounded if  $p_i(t)\in C^2,\, n+
1 \leq i \leq 2n;\, p_i(t) \in C^1, \,0 \leq i \leq n$. Recently, we \cite{Li} obtained the same conclusion if  $p_i(t)\in C^1,\, n+
1 \leq i \leq 2n;\, p_i(t) \in C^0, \,0 \leq i \leq n$.

There are other results about the boundedness problem for superlinear Duffing equations during the past years, see \cite{Levi,  Liu,  Liu92,  Wang2009,  WangWang2010} and the references therein.  As for constructing unbounded solutions for superlinear Duffing equations, there also are some results (\cite{Littlewood66,  Littlewood661,  Long,  Levi92,  Levi97, Wang2000}). Let us recall the results in \cite{Levi97} and \cite{Wang2000}. Levi and You \cite{Levi97} proved that the equation
$$
x'' + x^{2n+1} + p(t)x^{2m+1} = 0
$$
with a special discontinuous coefficient $p(t)=K^{[t]mod2}, \, 0 < K < 1, \, n+2\leq 2m + 1 <2n+1$, possesses an oscillatory unbounded solution. In 2000, Wang \cite{Wang2000}
constructed a continuous periodic function $p(t)$ such that the corresponding
equation
$$
x'' + x^{2n+1} + p(t)x^{i} = 0
$$
possesses a solution which escapes to infinity
in some finite time, where $n \geq 2$ and $n+2\leq i<2n+1$.

In this paper we consider the following  second order differential equation
\begin{equation}\label{w}
x''+a(x)\,x^{2n+1}+p(t)\,x^{2m+1}=0,\ \ \ \ n+2\leq 2m+1<2n+1,
\end{equation}
where the coefficient $a(x)$ is an arbitrary positive smooth periodic function defined in the whole real axis, will construct a continuous positive periodic function $p(t)$  and obtain the coexistence of quasi-periodic solutions and blow-up phenomena for the corresponding
equation \eqref{w}. More precisely, we will prove

\begin{theorem}\label{mainresult}
There exists a continuous positive periodic function $p(t)$ such that the corresponding
equation \eqref{w} possesses a solution which escapes to infinity in some finite time, and also has infinitely many subharmonic and quasi-periodic solutions.
\end{theorem}

Firstly, we will employ the idea in \cite{Wang2000} to construct a continuous positive periodic function $p(t)$ such that the corresponding
equation \eqref{w} possesses a solution which escapes to infinity in some finite time. Here, we will construct the positive periodic function $p(t)$ and the blow up solution $x(t)$ simultaneously.

First of all, we observe that during the time when the curve spirals once around the origin, the action variable $I$ increases at some times and decreases at other times after the action-angle variables $(I,\theta)$ are introduced. Therefore we do not know whether the increment of $I$ is positive or negative. However we can construct a time $t_{1}\ll1$ and modify $p(t)\equiv 1$ on $[0, 1]$ so that the increment of $I$ on this time interval $[0, t_{1}]$  is positive and equals to $O\left(\frac{1}{\tau}I_{0}^{\frac{2m-n+2}{n+2}}\right)$ if the  initial point $\big(I(0), \theta(0)\big)=(I_{0} , 0) $ is far enough from the origin,	where the "jump" $\frac{1}{\tau}\,(0<\frac{1}{\tau}< 1)$ is critical to modify $p(t)$ and to our estimations.

Inductively, we can construct a	series of times $t_{1}, t_{2}, \cdots, t_{i} , t_{i+1}, \cdots$ and modify $p(t)$ on $[t_{i} , t_{i+1}], i=1, 2, \cdot\cdot\cdot, \,$ so that on every such interval $[t_{i} , t_{i+1}],$  the increment is positive and at least $O\left(\frac{1}{\tau}I_{0}^{\frac{2m-n+2}{n+2}}\right).$

Hence, we can construct a time $T_{1} \leq\frac{1}{\tau'}<1,$ so that the curve spirals at least  $ \left[\frac{1}{\tau'}I_{0}^{\frac{n}{n+2}} \right]$ times around the origin on the interval $[0, T_{1}]$ and $I_{1}:=I(T_{1})>I_{0}+\frac{c}{\tau\tau'}I_{0}^{\frac{2m+2}{n+2}} $
 with $c>0$ independent of induction steps, where $\frac{2m+2}{n+2}>1$ and sufficiently large $\tau'$ is used to ensure the blow up time not more than $1$. This complete an induction step: during the interval of time $[0, T_{1}],\, I$  increases from $I_{0}$ to $I_{1}$.

Inductively, a series of times $T_{1}, T_{2}, \cdots, T_{i}, T_{i+1}, \cdots$ are constructed such that during the interval of time $[T_{k} , T_{k+1}], \,I$ increases from $I_{k}$ to
$I_{k+1},$ where $I_{k+1}>I_{k}+ \frac{c}{(\tau\tau')^{k}}I_{k}^{\frac{2m+2}{n+2}} $ with the jump $\frac{1}{\tau^{k}}, $
where $T_{k+1}-T_{k} \leq \frac{1}{\tau'^{k}}$. The reason that the jump is less and less is that we have to assure $p(t) $ is
continuous. Because the exponent $\frac{2m+2}{n+2}>1, $ the less and less jump cannot stop the rapid increase of $I$. If  $\frac{1}{\tau'}$ is chosen small enough, we will find that $ T_{k} \rightarrow T_{\infty} <1$ as $k \rightarrow \infty $ and  $ I(t) \rightarrow +\infty$ as $t \rightarrow T_{\infty}. $

Once we have found the continuous positive periodic function $p(t)$ such that the corresponding
equation \eqref{w} possesses a solution which escapes to infinity in some finite time, the remain thing is to apply the result in \cite{Liu96}. To this end, we first introduce this result. Consider the conservative system
\begin{equation}\label{liu}
x''+p(t)\,x^{2m+1}+e(t,x)=0,\ \ \ \ m\geq 1,
\end{equation}
where $p(t)$ is a continuous and $1$-periodic function in the time $t$, $e(t, x)$ is also
$1$-periodic in the time $t$ and dominated by the power $x^{2m+2}$ in a neighborhood of
$x=0$.
Liu in \cite{Liu96} proved that if $\int_0^1p(t)dt\neq0$, then the trivial solution $x=0$ of Eq.(\ref{liu}) is
stable in the Liapunov sense if and only if $\int_0^1p(t)dt>0$ by Moser's twist theorem. Moreover, by the same argument in \cite{Diekerhoff}, if $\int_0^1p(t)dt>0$, then Eq.(\ref{liu}) also has infinitely many subharmonic and quasi-periodic solutions with small amplitudes. Compared Eq.\eqref{w} with  Eq.(\ref{liu}), since $p(t)$ is positive, then Eq.\eqref{w} also has infinitely many subharmonic and quasi-periodic solutions with small amplitudes.

Therefore, if we find a continuous positive periodic function $p(t)$ such that the corresponding
equation \eqref{w} possesses a blow-up solution, then such equation also has infinitely many subharmonic and quasi-periodic solutions with small amplitudes simultaneously.

Similar to the above, if we modify $p(t)\equiv0$ in $[0,1]$, we can construct a  continuous non-positive periodic function $p(t)$ with $\int_0^1p(t)dt<0$ such that the corresponding
equation \eqref{w} possesses a solution which escapes to infinity in some finite time, and at the same time, by Liu's result in \cite{Liu96}, the trivial solution $x=0$ of such equation is not stable. Therefore we can obtain

\begin{theorem}\label{mainresult1}
There exists a continuous non-positive periodic function $p(t)$ such that the corresponding
equation \eqref{w} possesses a solution which escapes to infinity in some finite time, and the trivial solution $x=0$ is not stable.
\end{theorem}

Finally we remark that the authors in \cite{WangWang2012} also obtained the coexistence of quasi-periodic solutions and blow-up phenomena in a class of higher dimensional Duffing-type equations, and the author \cite{Maro} obtained the coexistence of  bounded and unbounded motions in a bouncing ball model.

The paper is organized as follows. The action-angle variables are introduced in Section 2. In Section 3 we first prove some Lemmas which will be useful later. After that, we will construct a continuous positive periodic function $p(t)$ and a series of times $T_{k}$, then obtain an unbounded solution of equation \eqref{w} and finish the proof of Theorem \ref{mainresult}.


\section{Action-angle variables}

In this section we first introduce action and angle variables. Let $y=x'$, then Eq.\eqref{w} is  equivalent to the following Hamiltonian system
\begin{equation*}
x'=\frac{\partial H}{\partial y},~~~~~~y'=-\frac{\partial H}{\partial x},
\end{equation*}
where the Hamiltonian is
\begin{equation}\label{Hl}
H(x,y,t)=\frac{1}{2}y^{2}+G(x)+\frac{p(t)}{2m+2}x^{2m+2}
\end{equation}
with
\begin{equation}\label{Gl}
G(x)=\int_{0}^{x}a(s)s^{2n+1}ds.
\end{equation}

In order to introduce action and angle variables, we consider the auxiliary autonomous system
$$x'=y,~~~~~~y'=-G'(x).$$
Since $a(x)>0$ for all $x\in\mathbb{R}$, then $G(x)>0$ for all $x\neq0$, and all solutions of this system are periodic. For every $h>0$, denote by $I(h)$ the area enclosed by the closed curve
$$
\Gamma_h:\ \ \dfrac{1}{2}y^{2}+G(x)=h.
$$
That is,
$$
I=I(h)=\int_{\Gamma_{h}}\sqrt{2(h-G(x))}dx.
$$
Let $h=h(I)$ be the  inverse function of $I=I(h).$ Define
$$S(x,I)=\left\{
\begin{array}{l}
\displaystyle\int_{x_{-}}^{x}\sqrt{2(h(I)-G(s))}ds, \ \ \ \ \ \ \ \,y\geq 0,\\[0.6cm]
I-\displaystyle\int_{x_{-}}^{x}\sqrt{2(h(I)-G(s))}ds,\ \ y<0,
\end{array}
\right.
$$
where $x_{-}=x_{-}(I)<0$ is determined uniquely by $G(x_{-})=h(I).$

Now we introduce the well-known action-angle transformation
$$y=\frac{\partial S}{\partial x},\ \ \ \ \ \ \theta=\frac{\partial S}{\partial I}.$$
Then
$$\theta=\left\{
\begin{array}{l}
h'(I)\displaystyle\int_{x_{-}}^{x}\frac{ds}{\sqrt{2(h(I)-G(s))}},\ \ \ \ \ \ \ \ \ \ y\geq 0,\\[0.6cm]
1-h'(I)\displaystyle\int_{x_{-}}^{x}\frac{ds}{\sqrt{2(h(I)-G(s))}},\ \ \ \ \ \,y<0.
\end{array}
\right.
$$
Denote
\begin{equation*}\label{f}
\Psi:\ (\theta,I)\rightarrow(x,y),
\end{equation*}
then under $\Psi$, the Hamiltonian $H$ of \eqref{Hl} is transformed into
\begin{equation*}
H_{1}=H\circ \Psi=h(I)+\frac{p(t)}{2m+2}x(I,\theta)^{2m+2},
\end{equation*}
and the corresponding Hamiltonian system is
\begin{equation}\label{dHl}
\left\{
\begin{array}{l}
\dfrac{d\theta}{dt}=h'(I)+p(t)x(I,\theta)^{2m+1}\partial_I x(I,\theta),\\[0.4cm]
\dfrac{dI}{dt}=-p(t)x(I,\theta)^{2m+1}\partial_\theta x(I,\theta).\end{array}
\right.
\end{equation}

In the following, we do not attempt to obtain estimates with particularly sharp constants.
Indeed, we suppress all constants, and use the notations $u \leq\cdot\ v$ and $u\ \cdot\leq v $ to indicate that
$u \leq c\, v$ and $c\, u \leq v$, respectively, with some constant $c > 0$.

Now we give some estimates on $h(I)$ and $x(I,\theta)$. For this purpose, we first give some properties of the potential function $G$.

\begin{lemma}\label{G}
For all $x\neq0$, we have
$$
 x^{2n+2}\ \cdot \leq G(x)\leq \cdot\ x^{2n+2},\ \ \ \ |G'(x)|\leq \cdot\ |x|^{2n+1},
$$
$$
|G''(x)|\leq \cdot\ \left(|x|^{2n+1}+x^{2n}\right), \ \ \left|\frac{G(x)}{G'(x)}\right|\leq\cdot\ |x|, \ \ \left|\frac{G(x)G''(x)}{G'(x)^2}\right|\leq\cdot\ (|x|+1).
$$
\end{lemma}
\Proof Since the periodic function $a(x)$ is positive, then there exist two positive constants $m, M$ such that $0< m\leq a(x)\leq M$ for all $x\in\mathbb{R}$ and from the expression (\ref{Gl}) of $G(x)$, we know that for $x>0$,
\begin{equation*}
\frac{m}{2n+2}x^{2n+2}\leq G(x)=\int_0^xa(s)s^{2n+1}ds\leq\frac{M}{2n+2}x^{2n+2};
\end{equation*}
and if $x<0$,
\begin{equation*}
\frac{m}{2n+2}x^{2n+2}\leq G(x)=\int_0^{-x}a(-s)s^{2n+1}ds\leq\frac{M}{2n+2}x^{2n+2},
\end{equation*}
which yields the first inequality of this lemma.

If $x>0$, then $mx^{2n+1}\leq G'(x)=a(x)x^{2n+1}\leq Mx^{2n+1}$; if $x<0$, then $Mx^{2n+1}\leq G'(x)=a(x)x^{2n+1}\leq mx^{2n+1}$. Combining this two inequalities, one can obtain the second inequality.

Since $G''(x)=a'(x)x^{2n+1}+(2n+1)a(x)x^{2n}$, if we let $m_1\leq a'(x)\leq M_1$ for all $x\in\mathbb{R}$, where $m_1<0, M_1>0$ are two constants, then $m_1x^{2n+1}+(2n+1)mx^{2n}\leq G''(x)\leq M_1x^{2n+1}+(2n+1)Mx^{2n}$ for $x>0$, and $M_1x^{2n+1}+(2n+1)mx^{2n}\leq G''(x)\leq m_1x^{2n+1}+(2n+1)Mx^{2n}$ for $x<0$, therefore $|G''(x)|\leq\max\{|m_1|,M_1\}|x|^{2n+1}+(2n+1)Mx^{2n}\leq\cdot\ \left(|x|^{2n+1}+x^{2n}\right)$, which is the third inequality.

For $x>0$, $\frac{m}{M}x\leq\frac{G(x)}{G'(x)}\leq\frac{M}{m}x$, and if $x<0$, $\frac{M}{m}x\leq\frac{G(x)}{G'(x)}\leq\frac{m}{M}x$, then the fourth inequality holds.

If $x>0$, $\frac{mm_1}{(2n+2)M^2}x+\frac{(2n+1)m^2}{(2n+2)M^2}\leq\frac{G(x)G''(x)}{G'(x)^2}\leq\frac{MM_1}{(2n+2)m^2}x+\frac{(2n+1)M^2}{(2n+2)m^2}$, and if  $x<0$, $\frac{mM_1}{(2n+2)M^2}x+\frac{(2n+1)m^2}{(2n+2)M^2}\leq\frac{G(x)G''(x)}{G'(x)^2}\leq\frac{Mm_1}{(2n+2)m^2}x+\frac{(2n+1)M^2}{(2n+2)m^2}$, hence
$$
\left|\frac{G(x)G''(x)}{G'(x)^2}\right|\leq\cdot\ (|x|+1).
$$
Up to now, we have finished the proof of the lemma. \qed

\begin{lemma}\label{I(h)}
For sufficiently large $h>0$, we have
$$
h^{\frac{1}{2}+\frac{1}{2n+2}}\ \cdot \leq I(h)\leq \cdot\ h^{\frac{1}{2}+\frac{1}{2n+2}},
$$
$$
h^{-\frac{1}{2}+\frac{1}{2n+2}}\ \cdot \leq I'(h)\leq\cdot\  h^{-\frac{1}{2}+\frac{1}{2n+2}},
$$
$$
|I''(h)|\leq\cdot\  h^{-\frac{3}{2}+\frac{1}{n+1}}.
$$
\end{lemma}
\Proof Let $x_{-}<0<x_+$ defined by $G(x_{-})=G(x_+)=h$ for any $h>0$, then
$$
I(h)=2\sqrt{2}\int_{x_{-}}^{x_+}\sqrt{h-G(x)}dx.
$$
By Lemma \ref{G}, we know that
$$
h^{\frac{1}{2n+2}}\ \cdot \leq x_+(h), |x_{-}(h)|\leq \cdot\ h^{\frac{1}{2n+2}},
$$
which implies that
$$
I(h)\leq \cdot\ h^{\frac{1}{2}+\frac{1}{2n+2}}.
$$
On the other hand, if we let $\bar{x}>0$ determined by $G(\bar{x})=\frac{h}{2}$, then
\begin{equation}\label{barx}
h^{\frac{1}{2n+2}}\ \cdot \leq\bar{x}(h)\leq \cdot\ h^{\frac{1}{2n+2}},
\end{equation}
and for $0\leq s\leq \bar{x}$,
$$
\sqrt{h-G(s)}\geq \sqrt{h-G(\bar{x})}\geq \cdot\ h^{\frac{1}{2}},
$$
thus
$$
I(h)\geq2\sqrt{2}\int_{0}^{\bar{x}}\sqrt{h-G(s)}ds\geq \cdot\ h^{\frac{1}{2}+\frac{1}{2n+2}}.
$$

Now we prove the second inequality. From the expression of $I(h)$, we know that
$$
I'(h)=\sqrt{2}\int_{x_{-}}^{x_+}\frac{dx}{\sqrt{h-G(x)}}=\sqrt{2}\int_{x_{-}}^{0}\frac{dx}{\sqrt{h-G(x)}}
+\sqrt{2}\int_{0}^{x_+}\frac{dx}{\sqrt{h-G(x)}}.
$$

The second term can be rewritten as follows
$$
\int_{0}^{x_+}\frac{dx}{\sqrt{h-G(x)}}=\int_{0}^{\bar{x}}\frac{dx}{\sqrt{h-G(x)}}+\int_{\bar{x}}^{x_+}\frac{dx}{\sqrt{h-G(x)}}.
$$
For $0\leq x\leq \bar{x}$, we have
$$
h^{-\frac{1}{2}}\ \cdot  \leq\frac{1}{\sqrt{h-G(x)}}\leq\cdot\ h^{-\frac{1}{2}},
$$
which together with (\ref{barx}) implies that
$$
h^{-\frac{1}{2}+\frac{1}{2n+2}}\ \cdot \leq\int_{0}^{\bar{x}}\frac{dx}{\sqrt{h-G(x)}}\leq\cdot\  h^{-\frac{1}{2}+\frac{1}{2n+2}}.
$$

If $\bar{x}\leq x\leq x_+$, then
$$
h-G(x)=G(x_+)-G(x)=G'(\xi)(x_+-x),\ \ \xi\in(x,x_+)\subset (\bar{x},x_+)
$$
and
$$
h^{1-\frac{1}{2n+2}}\ \cdot \leq \bar{x}^{2n+1}\ \cdot \leq G'(\xi)=a(\xi)\xi^{2n+1}\leq\cdot\ x_+^{2n+1}\leq\cdot\ h^{1-\frac{1}{2n+2}},
$$
which implies that
$$
\frac{h^{-\frac{1}{2}+\frac{1}{4n+4}}}{\sqrt{x_+-x}}\ \cdot  \leq\frac{1}{\sqrt{h-G(x)}}\leq\cdot\ \frac{h^{-\frac{1}{2}+\frac{1}{4n+4}}}{\sqrt{x_+-x}}.
$$
Combining the last equation with the fact that
$$
h^{\frac{1}{4n+4}}\ \cdot \leq\int_{\bar{x}}^{x_+}\frac{dx}{\sqrt{x_+-x}}\leq\cdot\  h^{\frac{1}{4n+4}},
$$
we obtain
$$
h^{-\frac{1}{2}+\frac{1}{2n+2}}\ \cdot \leq\int_{\bar{x}}^{x_+}\frac{dx}{\sqrt{h-G(x)}}\leq\cdot\  h^{-\frac{1}{2}+\frac{1}{2n+2}},
$$
and thus
$$
h^{-\frac{1}{2}+\frac{1}{2n+2}}\ \cdot \leq\int_{0}^{x_+}\frac{dx}{\sqrt{h-G(x)}}\leq\cdot\  h^{-\frac{1}{2}+\frac{1}{2n+2}}.
$$

Similarly, one can prove that
$$
h^{-\frac{1}{2}+\frac{1}{2n+2}}\ \cdot \leq\int_{x_{-}}^{0}\frac{dx}{\sqrt{h-G(x)}}\leq\cdot\  h^{-\frac{1}{2}+\frac{1}{2n+2}},
$$
which completes the proof of the second inequality.

Finally we prove the estimate on $I''(h)$. From \cite{Levi}, we know that
$$
I''(h)=\frac{\sqrt{2}}{h}\dint_{x_{-}}^{x_+}\left(\frac{1}{2}-\frac{G(x)G''(x)}{G'(x)^2}\right)\frac{dx}{\sqrt{h-G(x)}}.
$$
By Lemma \ref{G}, for sufficiently large $h>0$, we have
$$
\left|\frac{G(x)G''(x)}{G'(x)^2}\right|\leq\cdot \ (|x|+1)\leq\cdot\ h^{\frac{1}{2n+2}},
$$
and thus
$$
|I''(h)|\leq\cdot\ h^{\frac{1}{2n+2}-1}I'(h)\leq\cdot\  h^{-\frac{3}{2}+\frac{1}{n+1}},
$$
which finishes the proof of this lemma.\qed

Since $h=h(I)$ is the inverse function of $I=I(h)$, we immediately obtain
\begin{lemma}\label{h(I)}
For sufficiently large $I>0$, we have
$$
I^{\frac{2n+2}{n+2}}\ \cdot \leq h(I)\leq \cdot\ I^{\frac{2n+2}{n+2}},\ \
I^{\frac{n}{n+2}}\ \cdot \leq h'(I)\leq\cdot\  I^{\frac{n}{n+2}},\ \  |h''(I)|\leq\cdot\  I^{-\frac{1}{n+2}}.
$$
\end{lemma}

Now we give some estimates on $x(I,\theta)$.

\begin{lemma}\label{x(I)}
For sufficiently large $I>0$, we have
$$
\partial_\theta x> 0, \hbox {when}\  y> 0,\ \ \ \ \partial_\theta x<0, \hbox {when }\ y<0;
$$
$$
|x|\leq \cdot \ I^{\frac{1}{n+2}},\ \ \ \
\left|\partial_{\theta} x\right|\leq \cdot \ I^{\frac{1}{n+2}},\ \ \ \
\left|\partial_I x\right|\leq\cdot \ I^{-\frac{n}{n+2}}.
$$
\end{lemma}

\Proof The inequality $|x|\leq \cdot \ I^{\frac{1}{n+2}}$ is obvious. According to the definition of $\theta$, when $y>0$ (or $0<\theta< \frac{1}{2}$), we have
$$
\theta=h'(I)\displaystyle\int_{x_{-}}^{x}\frac{ds}{\sqrt{2(h(I)-G(s))}}.
$$
Differentiating the above equality with respect to $\theta$ yields that
$$
1=\frac{h'(I)}{\sqrt{2(h(I)-G(x))}}\partial_\theta x,
$$
which implies that
$$
\partial_\theta x=\frac{\sqrt{2(h(I)-G(x))}}{h'(I)}.
$$
Therefore, $\partial_\theta x>0$ for $0<\theta<\frac{1}{2}$, and by Lemma \ref{h(I)}, $\partial_{\theta} x\leq \cdot \ I^{\frac{1}{n+2}}$ holds.

Similarly, when $y<0$ (or $\frac{1}{2}<\theta<1$), we have
$$
\partial_\theta x=-\frac{\sqrt{2(h(I)-G(x))}}{h'(I)},
$$
hence $\partial_\theta x<0$ and $|\partial_{\theta} x|\leq \cdot \ I^{\frac{1}{n+2}}$.

Now we prove the estimate on $\partial_I x$. From \cite{Levi}, when $y\geq0$ (or $0\leq \theta\leq \frac{1}{2}$), we have
$$
\partial_I x=\sqrt{2(h(I)-G(x))}\displaystyle\int_{x_{-}}^{x}\frac{L(I,s) ds}{\sqrt{2(h(I)-G(s))}}+\frac{h'(I)}{h(I)}\frac{G(x)}{G'(x)},
$$
where
$$
L(I,x)=-\frac{h''(I)}{h'(I)}-\frac{h'(I)}{2h(I)}\left(1-\frac{2G(x)G''(x)}{G'(x)^2}\right).
$$
According to Lemmas \ref{G}, \ref{h(I)}, $|L|\leq\cdot \ I^{-\frac{n+1}{n+2}} $ and $\left|\frac{h'(I)}{h(I)}\frac{G(x)}{G'(x)}\right|\leq\cdot \ I^{-\frac{n+1}{n+2}}$ hold. Also, since for $x_{-}\leq s\leq x$, $\frac{\sqrt{2(h(I)-G(x))}}{\sqrt{2(h(I)-G(s))}}\leq 1$, then $|\partial_I x|\leq\cdot \ I^{-\frac{n}{n+2}}$. One can obtain the same estimate for $y<0$. Thus, we have finished the proof of this lemma.\qed

If we define $x_1, x_2, x_3$ by
$$
x(I,\theta)=I^{\frac{1}{n+2}}x_1(I,\theta),\ \  \partial_I x(I,\theta)=I^{-\frac{n}{n+2}}x_2(I,\theta), \ \ \partial_{\theta} x(I,\theta)=I^{\frac{1}{n+2}}x_3(I,\theta),
$$
then they are bounded functions for sufficiently large $I$ and all $\theta\in\mathbb{R}$, that is, there are three positive constants $B_1, B_2, B_3$ such that
\begin{equation}\label{B}
|x_i(I,\theta)|\leq B_i,\ \ \ \ i=1,2,3.
\end{equation}
Furthermore, there exist two positive constants $C_1, C_2$ such that
\begin{equation}\label{C}
-x_1(I, \theta)\geq C_1,\ \ \ \ x_3(I, \theta)\geq C_2
\end{equation}
for sufficiently large $I$ and $\theta\in[\frac{1}{16}, \frac{3}{16}]$. Moreover, we can rewrite system (\ref{dHl}) into
\begin{equation}\label{dH2}
\left\{
\begin{array}{l}
\dfrac{d\theta}{dt}=h'(I)+p(t)I^{\frac{2m+1-n}{n+2}}x_1(I,\theta)^{2m+1}x_2(I,\theta),\\[0.4cm]
\dfrac{dI}{dt}=-p(t)I^{\frac{2m+2}{n+2}}x_1(I,\theta)^{2m+1}x_3(I,\theta).\end{array}
\right.
\end{equation}
By Lemma \ref{h(I)}, we know that $I^{\frac{n}{n+2}}\ \cdot \leq h'(I)\leq\cdot\  I^{\frac{n}{n+2}}$. Also by our assumption, $2m+1\leq 2(n-1)+1$, thus
$2m+1-n\leq n-1$ and the right of the first equation in (\ref{dH2}) is dominated by $h'(I)$ for sufficiently large $I$.


\section{The proof of Theorem \ref{mainresult} }

Now we define $p(t)$ in $[0, 1]$.   We will construct a time $t_1<1$ and modify
$p^0(t)\equiv1$ on $[0, 1]$ so that the action of one solution of (\ref{dH2}) increases in $[0, t_1]$.
We divide the construction into two steps: first, we construct a piecewise
continuous function $p^1(t)$ so that the action of one solution of (\ref{dH2}) obtains a positive increment in $[0, t_1]$ as we expect. Then we
modify this function $p^1(t)$ into a continuous one in such a way that the modification does not influence the estimate we had obtained before.

Without loss of generality, we assume that the function $a(x)$ is even. Denote the corresponding Hamiltonian system (\ref{dH2}) with the coefficient function $p(t)$ by $X_p$. Suppose the solution $(I(t), \theta(t))$ of $X_{p^0}$ with $(I(0), \theta(0))=(I_0, 0)$ at $t=0$ arrives at $(I_{\frac{1}{4}}, \frac{1}{4})$ at $t=t_{\frac{1}{4}}\ll 1$, where $I_0$ is a sufficiently large constant which will be determined later. Define $p^{\frac{1}{4}}(t)$ be a piecewise continuous function as follows
\begin{equation*}
p^{\frac{1}{4}}(t)=\left\{
\begin{array}{lll}
&1,& t\in[0,t_{\frac{1}{4}}],\\[0.2cm]
&1-\sigma,&t\in(t_{\frac{1}{4}},1],
\end{array}
\right.
\end{equation*}
where $0<\sigma<1$ is the jump, which is used to
control the increment of $I$.

Suppose the solution $(I(t), \theta(t))$ of $X_{p^{\frac{1}{4}}}$ with $(I(0), \theta(0))=(I_0, 0)$ at $t=0$ arrives at $(I_{\frac{1}{2}}, \frac{1}{2})$ at $t=t_{\frac{1}{2}}\ll 1$.  Define $p^{\frac{1}{2}}(t)$ be a piecewise continuous function as follows
\begin{equation*}
p^{\frac{1}{2}}(t)=\left\{
\begin{array}{lll}
&1,& t\in[0,t_{\frac{1}{4}}],\\[0.2cm]
&1-\sigma,&t\in(t_{\frac{1}{4}},t_{\frac{1}{2}}],\\[0.2cm]
&1,& t\in(t_{\frac{1}{2}},1].
\end{array}
\right.
\end{equation*}

Suppose the solution $(I(t), \theta(t))$ of $X_{p^{\frac{1}{2}}}$ with $(I(0), \theta(0))=(I_0, 0)$ at $t=0$ arrives at $(I_{\frac{3}{4}}, \frac{3}{4})$ at $t=t_{\frac{3}{4}}\ll 1$.  Define $p^{\frac{3}{4}}(t)$ be a piecewise continuous function as follows
\begin{equation*}
p^{\frac{3}{4}}(t)=\left\{
\begin{array}{lll}
&1,& t\in[0,t_{\frac{1}{4}}],\\[0.2cm]
&1-\sigma,&t\in(t_{\frac{1}{4}},t_{\frac{1}{2}}],\\[0.2cm]
&1,& t\in(t_{\frac{1}{2}},t_{\frac{3}{4}}],\\[0.2cm]
&1,& t\in(t_{\frac{3}{4}}, 1].
\end{array}
\right.
\end{equation*}

Suppose the solution $(I(t), \theta(t))$ of $X_{p^{\frac{3}{4}}}$ with $(I(0), \theta(0))=(I_0, 0)$ at $t=0$ arrives at $(I_{1}, 1)$ at $t=t_{1}\ll 1$.  Define $p^1(t)$ be a piecewise continuous function as follows
\begin{equation}\label{pp}
p^{1}(t)=\left\{
\begin{array}{lll}
&1,& t\in[0,t_{\frac{1}{4}}],\\[0.2cm]
&1-\sigma,&t\in(t_{\frac{1}{4}},t_{\frac{1}{2}}],\\[0.2cm]
&1,& t\in(t_{\frac{1}{2}},t_{\frac{3}{4}}],\\[0.2cm]
&1-\sigma,& t\in(t_{\frac{3}{4}}, t_1],\\[0.2cm]
&1,& t\in (t_{1}, 1].
\end{array}
\right.
\end{equation}

That is to say, the solution $(I(t), \theta(t))$ of $X_{p^1}$ with $(I(0), \theta(0))=(I_0, 0)$ at $t=0$ arrives at $(I_{\frac{1}{4}}, \frac{1}{4})$ at $t=t_{\frac{1}{4}}$, arrives at $(I_{\frac{1}{2}}, \frac{1}{2})$ at $t=t_{\frac{1}{2}}$, arrives at $(I_{\frac{3}{4}}, \frac{3}{4})$ at $t=t_{\frac{3}{4}}$, arrives at $(I_{1}, 1)$ at $t=t_{1}$, which finishes one cycle of the construction of $p(t)$.

Now we estimate the differences $I_1-I_0$ and $t_1-t_0$.

\begin{lemma}\label{1/4}
If $I_{0}$ is sufficiently large, then
$$
I_{0}^{-\frac{n}{n+2}}\ \cdot  \leq t_{\frac{1}{4}}\leq\cdot \ I_{0}^{-\frac{n}{n+2}},$$
$$
I_{0}^{\frac{2m+2-n}{n+2}}\ \cdot  \leq I_{\frac{1}{4}}-I_{0}\leq\cdot \ I_{0}^{\frac{2m+2-n}{n+2}}.
$$
\end{lemma}
\Proof Because $x_1<0, x_3>0$ for $\theta\in(0,\frac{1}{4})$, then $\frac{dI}{dt}>0$ for $t\in(0,t_{\frac{1}{4}})$ and thus $I(t)$ is an increasing function in this interval. Integrating the first equation of (\ref{dH2}) from $t=0$ to $t=t_{\frac{1}{4}}$ yields that
\begin{equation}\label{t1/4}
t_{\frac{1}{4}}\leq \frac{1}{4}\left(h'(I_0)-B_1^{2m+1}B_2I_{0}^{\frac{2m+1-n}{n+2}}\right)^{-1}\leq\cdot \ I_{0}^{-\frac{n}{n+2}}
\end{equation}
for sufficiently large $I_0>0$, here we use the estimate on $h'(I)$ in Lemma \ref{h(I)} and the bound $B_i$ of $x_i$ in (\ref{B}).

From the second equation of (\ref{dH2}), we have
\begin{equation}\label{dI}
\dfrac{dI}{I^{\frac{2m+2}{n+2}}}=-x_1(I,\theta)^{2m+1}x_3(I,\theta)dt.
\end{equation}
Since $x_1<0, x_3>0$ for $\theta\in(0,\frac{1}{4})$, and also integrating the above equation (\ref{dI}) from $t=0$ to $t=t_{\frac{1}{4}}$ , one can obtain
\begin{equation}\label{dI0}
\frac{n+2}{2m-n}\left(I_{\frac{1}{4}}^{-\frac{2m-n}{n+2}}-I_{0}^{-\frac{2m-n}{n+2}}\right)=\int_0^{t_{\frac{1}{4}}}x_1(I,\theta)^{2m+1}x_3(I,\theta)dt,
\end{equation}
and by (\ref{t1/4}), we get
$$
\frac{n+2}{2m-n}\left(I_{0}^{-\frac{2m-n}{n+2}}-I_{\frac{1}{4}}^{-\frac{2m-n}{n+2}}\right)\leq cB_1^{2m+1}B_2I_0^{-\frac{n}{n+2}},
$$
where the constant $c>0$ is given by (\ref{t1/4}), which implies that
$$
I_{\frac{1}{4}}^{-\frac{2m-n}{n+2}}\geq I_{0}^{-\frac{2m-n}{n+2}}-\tilde{c}I_0^{-\frac{n}{n+2}}=I_{0}^{-\frac{2m-n}{n+2}}\left(1-\tilde{c}I_0^{-\frac{2n-2m}{n+2}}\right),
$$
here the constant $\tilde{c}=c\frac{2m-n}{ n+2}B_1^{2m+1}B_2>0$. Hence we obtain
$$
I_{\frac{1}{4}}\leq I_{0}\left(1-\tilde{c}I_0^{-\frac{2n-2m}{n+2}}\right)^{-\frac{n+2}{2m-n}},
$$
which leads to
$$
I_{\frac{1}{4}}-I_{0}\leq\cdot \ I_{0}^{\frac{2m+2-n}{n+2}}.
$$

On the other hand, from the first equation of (\ref{dH2}), we have
$$
\begin{array}{lll}
t_{\frac{1}{4}}&\geq& \frac{1}{4}\left(h'(I_{\frac{1}{4}})+B_1^{2m+1}B_2 I_{\frac{1}{4}}^{\frac{2m+1-n}{n+2}}\right)^{-1}\\[0.4cm]
               &\geq\ \cdot& \left(I_{\frac{1}{4}}^{\frac{n}{n+2}}+ I_{\frac{1}{4}}^{\frac{2m+1-n}{n+2}}\right)^{-1}\\[0.4cm]
                &\geq\ \cdot& I_{\frac{1}{4}}^{-\frac{n}{n+2}}\\[0.4cm]
                &\geq\ \cdot& I_{0}^{-\frac{n}{n+2}}.
\end{array}
$$

Finally, it follows from (\ref{C}) and (\ref{dI0}) that
$$
\frac{n+2}{2m-n}\left(I_{0}^{-\frac{2m-n}{n+2}}-I_{\frac{1}{4}}^{-\frac{2m-n}{n+2}}\right)\geq cC_1^{2m+1}C_2\left(t_{\frac{3}{16}}-t_{\frac{1}{16}}\right).
$$
Similarly, the following estimate
$$
t_{\frac{3}{16}}-t_{\frac{1}{16}}\geq \cdot\ I_{0}^{-\frac{n}{n+2}}
$$
holds. Combining the two inequalities above yields that
$$
I_{0}^{-\frac{2m-n}{n+2}}-I_{\frac{1}{4}}^{-\frac{2m-n}{n+2}}\geq \cdot\ I_{0}^{-\frac{n}{n+2}},
$$
and
$$
I_{\frac{1}{4}}-I_{0}\geq\cdot \ I_{0}^{\frac{2m+2-n}{n+2}}.
$$
\qed

\begin{lemma}\label{1/2}
If $I_{0}$ is sufficiently large, then
$$
I_{0}^{-\frac{n}{n+2}}\ \cdot  \leq t_{\frac{1}{2}}- t_{\frac{1}{4}}\leq\cdot \ I_{0}^{-\frac{n}{n+2}},$$
$$
(1-\sigma)I_{0}^{\frac{2m+2-n}{n+2}}\ \cdot  \leq I_{\frac{1}{4}}-I_{\frac{1}{2}}\leq\cdot \ (1-\sigma)I_{0}^{\frac{2m+2-n}{n+2}}.
$$
\end{lemma}
\Proof Because $x_1>0, x_3>0$ for $\theta\in(\frac{1}{4},\frac{1}{2})$, then $\frac{dI}{dt}<0$ for $t\in(t_{\frac{1}{4}},t_{\frac{1}{2}})$ and thus $I(t)$ is an decreasing function in this interval. Integrating the first equation of (\ref{dH2}) from $t=t_{\frac{1}{4}}$ to $t=t_{\frac{1}{2}}$ yields that
\begin{equation}\label{t01/4}
I_{\frac{1}{4}}^{-\frac{n}{n+2}}\ \cdot  \leq t_{\frac{1}{2}}- t_{\frac{1}{4}}\leq\cdot \ I_{\frac{1}{2}}^{-\frac{n}{n+2}}.
\end{equation}

From the second equation of (\ref{dH2}), we have
\begin{equation*}
(1-\sigma) I_{\frac{1}{4}}^{-\frac{n}{n+2}} \ \cdot \leq I_{\frac{1}{2}}^{-\frac{2m-n}{n+2}}-I_{\frac{1}{4}}^{-\frac{2m-n}{n+2}}\leq\cdot \ (1-\sigma) I_{\frac{1}{2}}^{-\frac{n}{n+2}},
\end{equation*}
and thus
\begin{equation}\label{t21/4}
(1-\sigma)I_{\frac{1}{4}}^{\frac{2m+2-n}{n+2}}\ \cdot  \leq I_{\frac{1}{4}}-I_{\frac{1}{2}}\leq\cdot \ (1-\sigma)I_{\frac{1}{2}}^{\frac{2m+2-n}{n+2}}.
\end{equation}

By Lemma \ref{1/4}, we have
$$
I_{\frac{1}{4}}\leq I_0\left(1+cI_0^{\frac{2m-2n}{n+2}}\right)
$$
with some constant $c>0$, which implies that
\begin{equation}\label{t31/4}
I_{\frac{1}{4}}^{-\frac{n}{n+2}}\geq \cdot\ I_{0}^{-\frac{n}{n+2}}
\end{equation}
holds for sufficiently large $I_0>0$. Meanwhile it follows from (\ref{t21/4}) that
$$
I_{\frac{1}{2}}\geq \cdot\ I_{\frac{1}{4}},
$$
which together with $I_{\frac{1}{4}}>I_0$ implies that
\begin{equation}\label{t51/4}
I_{\frac{1}{2}}\geq \cdot\ I_{0},
\end{equation}
and
\begin{equation}\label{t41/4}
I_{\frac{1}{2}}^{-\frac{n}{n+2}}\leq \cdot\ I_{0}^{-\frac{n}{n+2}}.
\end{equation}
Combining (\ref{t01/4}), (\ref{t31/4}) with  (\ref{t41/4}), we obtain
$$
I_{0}^{-\frac{n}{n+2}}\ \cdot  \leq t_{\frac{1}{2}}- t_{\frac{1}{4}}\leq\cdot \ I_{0}^{-\frac{n}{n+2}}.
$$
Also, according to (\ref{t21/4}), (\ref{t51/4}) and $I_{\frac{1}{4}}\geq I_0$, one can obtain the second inequality in this lemma. \qed

Using the same method, one can prove the following result.

\begin{lemma}\label{1/1}
If $I_{0}$ is sufficiently large, then
$$
I_{0}^{-\frac{n}{n+2}}\ \cdot  \leq t_{\frac{3}{4}}- t_{\frac{1}{2}}\leq\cdot \ I_{0}^{-\frac{n}{n+2}},
$$
$$
I_{0}^{-\frac{n}{n+2}}\ \cdot  \leq t_{1}- t_{\frac{3}{4}}\leq\cdot \ I_{0}^{-\frac{n}{n+2}},
$$
$$
I_{0}^{\frac{2m+2-n}{n+2}}\ \cdot  \leq I_{\frac{3}{4}}-I_{\frac{1}{2}}\leq\cdot \ I_{0}^{\frac{2m+2-n}{n+2}},
$$
$$
(1-\sigma)I_{0}^{\frac{2m+2-n}{n+2}}\ \cdot  \leq I_{\frac{3}{4}}-I_{1}\leq\cdot \ (1-\sigma)I_{0}^{\frac{2m+2-n}{n+2}}.
$$
\end{lemma}

Combining Lemmas \ref{1/4}, \ref{1/2} and \ref{1/1}, we can obtain immediately the estimates on the time $t_1$ when the
curve spirals once around the origin and  the increment of the action variable $I_1-I_0$.

\begin{lemma}\label{tI}
If $I_{0}$ is sufficiently large, then
$$
I_{0}^{-\frac{n}{n+2}}\ \cdot  \leq t_{1}\leq\cdot \ I_{0}^{-\frac{n}{n+2}},$$
$$
\sigma I_{0}^{\frac{2m+2-n}{n+2}}\ \cdot  \leq I_{1}-I_{0}\leq\cdot \ \sigma I_{0}^{\frac{2m+2-n}{n+2}}.
$$
\end{lemma}

Now we modify the piecewise continuous function $p^1(t)$ of \eqref{pp} into a continuous one. Being short of signs, we keep the notations unchanged in the process of modification. For example, $p^1(t)$ denotes the continuous function modified from the original piecewise continuous function $p^1(t)$.

First we modify  $p^1(t)$ on the interval $[t_{\frac{1}{4}}, t_{\frac{1}{4}}+I_0^{-\eta}]$ $(\eta>\frac{n}{n+2})$ to be $\sigma(t_{\frac{1}{4}}-t)I_0^{\eta}+1$. It is easy to see that $\{(t, p^1(t)): t\in[t_{\frac{1}{4}}, t_{\frac{1}{4}}+I_0^{-\eta}]\}$ is the line segment connecting $(t_{\frac{1}{4}} , 1)$ and $(t_{\frac{1}{4}}+I_0^{-\eta},1-\sigma)$.

In view of the mean value theorem, there must exist a unique new time $t_{\frac{1}{2}}$ such that $\theta(t_{\frac{1}{2}})=\frac{1}{2}$  if we let
\begin{equation*}
p^{\frac{1}{2}}(t)=\left\{
\begin{array}{lll}
&1,& t\in[0,t_{\frac{1}{4}}],\\[0.2cm]
&\sigma(t_{\frac{1}{4}}-t)I_0^{\eta}+1,& t\in(t_{\frac{1}{4}}, t_{\frac{1}{4}}+I_0^{-\eta}],\\[0.2cm]
&1-\sigma,&t\in(t_{\frac{1}{4}}+I_0^{-\eta},t_{\frac{1}{2}}-I_0^{-\eta}],\\[0.2cm]
&\sigma(t-t_{\frac{1}{2}})I_0^{\eta}+1,& t\in(t_{\frac{1}{2}}-I_0^{-\eta}, t_{\frac{1}{2}}],\\[0.2cm]
&1,& t\in(t_{\frac{1}{2}},1].
\end{array}
\right.
\end{equation*}

Similarly, there exist the unique $t_{\frac{3}{4}}$ and $t_1$ such that $\theta(t_{\frac{3}{4}})=\frac{3}{4}$ and $\theta(t_1)=1$ for $X_{p^1}$ with
\begin{equation*}
p^{1}(t)=\left\{
\begin{array}{lll}
&1,& t\in[0,t_{\frac{1}{4}}],\\[0.2cm]
&\sigma(t_{\frac{1}{4}}-t)I_0^{\eta}+1,& t\in(t_{\frac{1}{4}}, t_{\frac{1}{4}}+I_0^{-\eta}],\\[0.2cm]
&1-\sigma,&t\in(t_{\frac{1}{4}}+I_0^{-\eta},t_{\frac{1}{2}}-I_0^{-\eta}],\\[0.2cm]
&\sigma(t-t_{\frac{1}{2}})I_0^{\eta}+1,& t\in(t_{\frac{1}{2}}-I_0^{-\eta}, t_{\frac{1}{2}}],\\[0.2cm]
&1,& t\in(t_{\frac{1}{2}},t_{\frac{3}{4}}],\\[0.2cm]
&\sigma(t_{\frac{3}{4}}-t)I_0^{\eta}+1,& t\in(t_{\frac{3}{4}}, t_{\frac{3}{4}}+I_0^{-\eta}],\\[0.2cm]
&1-\sigma,&t\in(t_{\frac{3}{4}}+I_0^{-\eta},t_{1}-I_0^{-\eta}],\\[0.2cm]
&\sigma(t-t_{1})I_0^{\eta}+1,& t\in(t_{1}-I_0^{-\eta}, t_{1}],\\[0.2cm]
&1, &t\in(t_{1}, 1].
\end{array}
\right.
\end{equation*}
Now the newest coefficient is already a continuous function. It is easy to check that Lemmas \ref{1/4}--\ref{tI} still hold with different constants after this modification in view of $I_{0}^{-\eta}\ll I_{0}^{-\frac{n}{n+2}}.$

We will modify $p^0$ inductively and denote the function obtained and the
corresponding solution with $(I_0, 0)$ as the initial point by $p^i$ and $(I^i(t), \theta^i(t))$ with $(I^i(t_i), \theta^i(t_i))=(I_i, i)$, respectively.

Suppose that we have obtained $p^0, p^1, \cdots, p^i$. The function $p^{i+1}$ defined on $[0,1]$ is constructed by modifying $p^i$ on the interval $[t_i, t_{i+1}]$, where $t_{i+1}$ satisfies $\theta^{i+1}(t_{i+1})=i+1$ in the same way as above if we regard $I_i, t_i$  as $I_0, t_0$. All the lemmas are true after the modification.

In the process of constructing $p^i$, we keep the jump $\sigma=\frac{1}{\tau}\, (\tau\geq2)$
unchanged until $i=j_1$. Then we let $\sigma=\frac{1}{\tau^2}$ and keep it unchanged until
$i=j_2$. Inductively, we choose $\sigma=\frac{1}{\tau^k}$ when $\theta\in [j_{k-1}, j_k]$, where $j_0=0,\,
j_1, j_2, \cdots$ are defined as below.

Let $j_{1}=\left[\frac{1}{\tau'}I_{0}^{\frac{n}{n+2}}\right]$, where $[x]$ denotes the integer part of $x$  and $\tau'>0$ is used to
control time and will be determined later. It follows that
$$T_{1}:=t_{j_{1}}\leq j_{1}\cdot t_{1} \ \cdot\leq
\frac{1}{\tau'}I_{0}^{\frac{n}{n+2}}\cdot I_{0}^{-\frac{n}{n+2}}
 \ \cdot \leq\frac{1}{\tau'}.$$
On the interval $[0, T_1]$, since $I_{k+1}-I_k\geq\cdot \ \frac{1}{\tau} I_{0}^{\frac{2m+2-n}{n+2}}$ for $k=0, 1, \cdots, j_1-1$, we have
$$
I_{j_1}:=I^{j_1}(T_1)\geq\cdot \ \left[\frac{1}{\tau'}I_{0}^{\frac{n}{n+2}}\right]\cdot\frac{1}{\tau} I_{0}^{\frac{2m+2-n}{n+2}}\geq\cdot \ \frac{1}{\tau\tau'}I_{0}^{\frac{2m+2}{n+2}}.
$$

Suppose we have defined $T_0=0, T_1, \cdots, T_i$ according to the above method and the
following are tenable for $ k=1,2,\cdots,i $:
$$
j_{k}-j_{k-1}=
\left[\frac{1}{\tau'^{k}}I_{j_{k-1}}^{\frac{n}{n+2}}\right], \ \ T_k-T_{k-1}\ \cdot \leq\frac{1}{\tau'^k},
$$
$$
\sigma=\tau_k=\frac{1}{\tau^{k}}, \ \ I_{j_k}:=I^{j_k}(T_k)\geq\cdot \ \frac{1}{(\tau\tau')^k}I_{0}^{\frac{2m+2}{n+2}}.
$$

Set
$$
j_{i+1}-j_{i}=
\left[\frac{1}{\tau'^{i+1}}I_{j_{i}}^{\frac{n}{n+2}}\right],\ \ \sigma=\tau_{i+1}=\frac{1}{\tau^{i+1}},\ \ T_{i+1}:=t_{j_{i+1}},
$$
similar to the above discussion, we have
$$
 I_{j_{i+1}}:=I^{j_{i+1}}(T_{i+1})\geq\cdot \ \frac{1}{(\tau\tau')^{i+1}}I_{0}^{\frac{2m+2}{n+2}},\ \ T_{i+1}-T_{i}\ \cdot \leq\frac{1}{\tau'^{i+1}}.
$$
Consequently,
$$
T_{i+1}\ \cdot \leq\sum_{k=1}^{i+1}\frac{1}{\tau'^{k}}\ \cdot \leq\frac{1}{\tau'}<1
$$
if $\tau'>0$ is sufficiently large.

Let
$$
\lim_{i\to\infty}T_i=T_{\infty}, \ \ \lim_{i\to\infty}p^i(t)=p(t),
$$
since
$$
\max_{t_1,t_2\in[T_k,T_{\infty}]}|p(t_1)-p(t_2)|\leq \frac{1}{\tau^{k}}, \ \ \lim_{t\to t_{\infty}}p(t)=1,
$$
then $p(t)$ can be extended to a continuous positive $1$-periodic function.

\begin{lemma}\label{4.2} If $I_{0}$ is sufficiently large, then
$$I_{j_{k}}\geq\cdot \ I_{0}^{l^{k}},$$
where the constant $l>1$.
\end{lemma}

\Proof  First, by the assumption $2m+1\ge n+2$,  if we let $l=\frac{2m+1}{n+2}+\frac{1}{2(n+2)}$, then $l>1$, and  for sufficiently large $I_{0}$ we have
$$I_{j_{1}}\geq\cdot \ \frac{1}{\tau\tau'}I_{0}^{\frac{2m+2}{n+2}}\geq\cdot \ \frac{1}{\tau\tau'}I_{0}^{\frac{1}{2(n+2)}}I_{0}^{l}\geq\cdot \ I_{0}^{l}.
$$
If
$$
I_{j_{k}}\geq\cdot \ I_{0}^{l^{k}},
$$
then
$$\begin{array}{lll}
I_{j_{k+1}}&\geq\cdot \ &\dfrac{1}{(\tau\tau')^{k+1}}I_{j_{k}}^{\frac{2m+2}{n+2}}\\[0.4cm]
&\geq\cdot \ & \dfrac{1}{(\tau\tau')^{k+1}}I_{0}^{\frac{2m+2}{n+2}l^{k}}\\[0.4cm]
&\geq\cdot \ & \dfrac{I_{0}^{\frac{l^{k}}{2(n+2)}}}{(\tau\tau')^{k+1}}I_{0}^{l^{k+1}}\\[0.5cm]
&\geq\cdot \ & I_{0}^{l^{k+1}}.
\end{array}$$
\qed

\noindent {\bf Proof of Theorem \ref{mainresult}}
\quad By Lemma \ref{4.2},  one has
$$\min_{t\in [T_{i},T_{i+1}]}I(t)\geq\cdot \ I_{j_{i}}\geq\cdot \ I_{0}^{l^{i}}.$$

Since $l>1$, then $I(t)\rightarrow +\infty$ as $t \rightarrow t_{\infty}$. Therefore,
Eq.\eqref{w} in Theorem \ref{mainresult} possesses an unbounded solution defined in the interval $[0, T_{\infty})$, and Theorem \ref{mainresult} is proved. \qed


\end{document}